\documentclass[12pt]{article}
\usepackage{amsmath,amsfonts,latexsym,amsthm,amssymb}
\newcommand{\F}{\mathbf F_q}
\newcommand{\K}{\overline{K}_c}
\newcommand{\A}{\mathcal A}
\newcommand{\G}{\mbox{gr}\, (\mathcal A)}
\newcommand{\LL}{\mathcal L}
\newcommand{\GL}{\mbox{gr}\, (\mathcal L)}
\begin{document}
\newtheorem*{lem}{Lemma}
\newtheorem{teo}{Theorem}
\pagestyle{plain}
\title{Differential Equations for $\F$-Linear Functions}
\author{Anatoly N. Kochubei\\ \footnotesize Institute of Mathematics,\\
\footnotesize Ukrainian National Academy of Sciences,\\
\footnotesize Tereshchenkivska 3, Kiev, 01601 Ukraine
\\ \footnotesize E-mail: \ ank@ank.kiev.ua}
\date{}
\maketitle
\vspace*{2cm}
Copyright Notice: This work has been submitted to Academic Press
for possible publication. Copyright may be transferred without
notice, after which this version may no longer be accessible.
\newpage
\vspace*{8cm}
\begin{abstract}
We study certain classes of equations for $\F$-linear functions, which
are the natural function field counterparts of linear ordinary
differential equations. It is shown that, in contrast to both
classical and $p$-adic cases, formal power series solutions have
positive radii of convergence near a singular point of an
equation. Algebraic properties of the ring of $\F$-linear
differential operators are also studied.
\end{abstract}
\vspace{2cm}
\newpage
\section{INTRODUCTION}
A standard model of a non-discrete, locally compact field of
a positive characteristic $p$ is the field $K$ of formal
Laurent series with coefficients from the Galois field $\F$,
$q=p^\upsilon $, $\upsilon \in \mathbf Z_+$. Let $\K$ be a
completion of an algebraic closure of $K$.

A function defined on a $\F$-subspace $K_0$ of $K$, with values
in $\K$, is called $\F$-linear
if $f(t_1+t_2)=f(t_1)+f(t_2)$ and $f(\alpha t)=\alpha f(t)$ for any
$t,t_1,t_2\in K$, $\alpha \in \F$.

$\F$-linear functions often appear in analysis over $\K$; see, in
particular, the works by Carlitz \cite{C1,C2}, Wagner \cite{W},
Goss \cite{G1,G2}, Thakur \cite{T1,T2}, and the author \cite{K1,K2}.
This class of functions includes, in particular, analogues of the
exponential, logarithm, Bessel, and hypergeometric functions.

It has been noticed (see e.g. \cite{T2,K1}) that the above functions
satisfy some equations, which can be seen as function field analogues of
first and second order linear differential equations with polynomial
coefficients. The role of a derivative is played by the operator
$$
d=\sqrt[q]{\ }\circ \Delta
$$
where $(\Delta u)(t)=u(xt)-xu(t)$, $x$ is a prime element in $K$.
The operator $d$ is also basic in the calculus of $\F$-linear functions
developed in \cite{K2}.

The meaning of a polynomial coefficient in the function field
case is not a usual multiplication by a polynomial, but the
action of a polynomial in the operator $\tau$, $\tau u=u^q$. This
was most clearly demonstrated in \cite{T2}, where an equation for
the hypergeometric functions was derived.

The aim of this paper is to start a general theory of
differential equations of the above kind. In fact we consider
equations (or systems) with holomorphic coefficients.
As in the classical theory, we have to make a distinction between
the regular and singular cases.

In the analytic theory of linear differential equations over
$\mathbb C$ a regular equation has a constant leading coefficient
(which can be assumed equal to 1). A leading coefficient of a
singular equation is a holomorphic function having zeroes at some
points. One can divide the equation by its leading coefficient,
but then poles would appear at other coefficients, and the
solution can have singularities (not only poles but in general
also essential singularities) at those points.

Similarly, in our case we understand a regular equation as the
one with the coefficient 1 at the highest order derivative.
As usual, a regular higher-order equation can be transformed
into a regular first-order system. For the regular case we obtain
a local existence and uniqueness theorem,
which is similar to analogous results for equations over $\mathbb
C$ or $Q_p$ (for the latter see \cite{Lu}). The only difference
is a formulation of the initial condition, which is specific for
the function field case.

The leading coefficient $A_m(\tau )$ of a singular $\F$-linear
equation of an order $m$ is a non-constant holomorhic function
of the operator $\tau$. Now one cannot divide the equation
$$
A_m(\tau )d^mu(t)+A_{m-1}(\tau )d^{m-1}u(t)+\ldots +A_0(\tau
)u(t)=f(t)
$$
for an $\F$-linear function $u(t)$ (note that automatically
$u(0)=0$) by $A_m(\tau )$. If $A_m(\tau )=\sum\limits_{i=0}^\infty
a_{mi}\tau ^i$, $a_{mi}\in \K$, then
$$
A_m(\tau )d^mu=
\sum\limits_{i=0}^\infty a_{mi}\left( d^mu\right) ^{q^i},
$$
and even when $A_m$ is a polynomial, in order to resolve our equation
with respect to $d^mu$ one has to solve an algebraic equation.

Thus for the singular case the situation looks even more
complicated than in the classical theory. However we show that
the behavior of the solutions cannot be too intricate.
Namely, in a striking contrast to the classical theory, any formal series
solution converges in some (sufficiently small) neighbourhood of
the singular point $t=0$. Note that in the $p$-adic case a similar
phenomenon takes place for equations satisfying certain strong conditions
upon zeros of indicial polynomials \cite{Cla,Set,Bald,Put1}. In our case
such a behavior is proved for any equation, which resembles the
(much simpler) case \cite{Put1} of differential equations over a
field of characteristics zero, whose residue field also has
characteristic zero.

We also study some algebraic properties of the ring of all
polynomial differential operators, that is the ring generated by
$\K$, $\tau$, and $d$. It is interesting to compare our results
with the ones for the case of characteristic zero \cite{Bj,D},
and the ones for usual differential operators over a field of
positive characteristic \cite{Put2}. It appears that the main
difference from the former case is caused by nonlinearity of
$\tau$ and $d$, while the latter case is totally different. For
example, the centre of our ring equals $\F$, the centre of the
ring of polynomial differential operators over a field $k$ with
$\mbox{char }k=0$ equals $k$. Meanwhile for the situation studied
in \cite{Put2} the centre is a ``big'' polynomial ring.

The author is grateful to David Goss for his constructive
criticism, which helped much to improve the exposition.

\section{ANALYTIC PROPERTIES}

Let us introduce some notation. If $t\in K$,
$$
t=\sum _{i=n}^\infty \theta_ix^i,\quad n\in \mathbf Z,\ \theta_i\in \F \ ,
\ \theta_n\ne 0,
$$
the absolute value $|t|$ is defined as
$$
|t|=q^{-n}.
$$
We preserve the notation $|\cdot |$ for the extension of the
absolute value onto $\K$. The norm $|P|$ of a matrix $P$ with
elements from $\K$ is defined as the maximum of absolute values
of the elements.

We will use systematically the Carlitz factorial $D_i$, $i\ge 0$,
defined as
$$
D_i=[i][i-1]^q\ldots [1]^{q^{i-1}},\ i\ge 1;\ D_0=1,
$$
where $[i]=x^{q^i}-x$ (this notation should not be confused with the one
for the commutator $[\cdot ,\cdot ]$). It is easy to see \cite{K2} that
\begin{equation}
d\left( \frac{t^{q^i}}{D_i}\right) = \frac{t^{q^{i-1}}}{D_{i-1}},\ i\ge
1;\quad d(\mbox{const})=0;
\end{equation}
\begin{equation}
\tau \left( \frac{t^{q^{i-1}}}{D_{i-1}}\right) =[i]\frac{t^{q^i}}{D_i},\ i\ge 1.
\end{equation}
It is known \cite{K2,T2} that $[d,\tau ]=[1]^{1/q}$.

\subsection{{\it Equations without Singularities}}

Let us consider an equation
\begin{equation}
dy(t)=P(\tau )y(t)+f(t)
\end{equation}
where for each $z\in \left(\K\right) ^m$, $t\in K$,
\begin{equation}
P(\tau )z=\sum \limits _{k=0}^\infty \pi _kz^{q^k},\quad
f(t)=\sum \limits _{j=0}^\infty \varphi _j\frac{t^{q^j}}{D_j},
\end{equation}
$\pi _k$ are $m\times m$ matrices with elements from $\K$,
$\varphi _j\in \left(\K\right) ^m$, and it is assumed that the series (4)
have positive radii of convergence. The action of the operator $\tau $
upon a vector or a matrix is defined component-wise, so that
$z^{q^k}=\left( z_1^{q^k},\ldots ,z_m^{q^k}\right)$ for
$z=(z_1,\ldots ,z_m)$. Similarly, if $\pi =(\pi _{ij})$ is a
matrix, we write $\pi ^{q^k}=\left( \pi_{ij}^{q^k}\right)$.

We will seek a $\F$-linear solution of (3) in some neighbourhood of the
origin, of the form
\begin{equation}
y(t)=\sum \limits _{i=0}^\infty y_i\frac{t^{q^i}}{D_i},\quad
y_i\in \left(\K\right) ^m,
\end{equation}
where $y_0$ is a given element, so that the ``initial'' condition
for our situation is
\begin{equation}
\lim \limits _{t\to 0}t^{-1}y(t)=y_0.
\end{equation}
Note that a function (5), provided the series has a positive
radius of convergence, tends to zero for $t\to 0$, so that the
right-hand side of (3) makes sense for small $|t|$.

\begin{teo}
For any $y_0\in \left(\K\right) ^m$ the equation (3) has a unique
local solution of the form (5), which satisfies (6), with the
series having a positive radius of convergence.
\end{teo}

{\it Proof}. Making (if necessary) the substitutions $t=c_1t'$,
$y=c_2y'$, with sufficiently small $|c_1|$, $|c_2|$, we may assume that
the coefficients in (4) are such that $\varphi _j\to 0$ for $j\to \infty$,
\begin{equation}
|\pi _k^q|\cdot q^{-\frac{q^{k+1}-q}{q-1}}\le 1,\quad k=0,1,\ldots .
\end{equation}

Using (1) and (2) we substitute (5) into (3), which results in the
recurrent formula for the coefficients $y_i$:
\begin{equation}
y_{l+1}=\sum \limits _{n+k=l}\pi
_k^qy_n^{q^{k+1}}[n+1]^{q^k}\ldots [n+k]^q+\varphi _l^q, \quad
l=0,1,2,\ldots ,
\end{equation}
where the expressions in square brackets are omitted if $k=0$.

It is seen from (8) that a solution of (3), (6) (if it exists)
is unique. Since $|[n]|=q^{-1}$ for all $n>0$, we find that
$$
\left| [n+1]^{q^k}\ldots [n+k]^q\right| =q^{-(q^k+\cdots +q)}=
q^{-\frac{q^{k+1}-q}{q-1}},
$$
and it follows from (7),(8) that
$$
|y_{l+1}|\le \max \left\{ |\varphi
_l|^q,|y_0|^{q^{l+1}},|y_1|^{q^l},\ldots ,|y_l|^q\right\} .
$$

Since $\varphi _n\to 0$, there exists such a number $l_0$ that
$|\varphi _l|\le 1$ for $l\ge l_0$. Now either $|y_l|\le 1$ for
all $l\ge l_0$ (and then the series (5) is convergent in a
neighbourhood of the origin), or $|y_{l_1}|>1$ for some $l_1\ge
l_0$. In the latter case
$$
|y_{l+1}|\le \max \left\{ |y_0|^{q^{l+1}},|y_1|^{q^l},\ldots ,
|y_l|^q\right\} ,\quad l\ge l_1.
$$

Let us choose $A>0$ in such a way that
$$
|y_l|\le A^{q^l},\quad l=1,2,\ldots ,l_1.
$$
Then it follows easily by induction that $|y_l|\le A^{q^l}$ for
all $l$, which implies the convergence of (5) near the
origin.$\quad \blacksquare $

\subsection{{\it Singular Equations}}

We will consider scalar equations of arbitrary order
\begin{equation}
\sum \limits _{j=0}^mA_j(\tau )d^ju=f
\end{equation}
where
$$
f(t)=\sum \limits _{n=0}^\infty \varphi _n\frac{t^{q^n}}{D_n},
$$
$A_j(\tau )$ are power series having (as well as the one for $f$)
positive radii of convergence.

It will be convenient to start from the model equation
\begin{equation}
\sum \limits _{j=0}^ma_j\tau ^jd^ju=f,\quad a_j\in \K,\ a_m\ne 0.
\end{equation}
Suppose that $u(t)$ is a formal solution of (10), of the form
\begin{equation}
u(t)=\sum \limits _{n=0}^\infty u_n\frac{t^{q^n}}{D_n}.
\end{equation}
Then
$$
a_0\sum \limits _{n=0}^\infty u_n\frac{t^{q^n}}{D_n}+\sum \limits
_{j=1}^ma_j\sum \limits _{n=j}^\infty u_n[n-j+1]\ldots [n]
\frac{t^{q^n}}{D_n}=\sum \limits _{n=0}^\infty \varphi
_n\frac{t^{q^n}}{D_n}.
$$
Changing the order of summation we find that for $n\ge m$
\begin{equation}
u_n\left( a_0+\sum \limits _{j=1}^ma_j[n-j+1]\ldots [n]\right) =
\varphi _n.
\end{equation}

Let us consider the expression
$$
\Phi _n=a_0+\sum \limits _{j=1}^ma_j[n-j+1]\ldots [n], \quad n\ge
m.
$$
Using repeatedly the identity $[i]^q+[1]=[i+1]$ we find that
$$
\Phi _n^{q^m}=a_o^{q^m}+\sum \limits _{j=1}^ma_j^{q^m}\prod
\limits _{k=0}^{j-1}[n-k]^{q^m}=
a_o^{q^m}+\sum \limits _{j=1}^ma_j^{q^m}\prod
\limits _{k=0}^{j-1}\left( [n]^{q^{m-k}}-\sum \limits
_{l=1}^k[1]^{q^{m-l}}\right),
$$
that is $\Phi _n^{q^m}=\Phi ^{(m)}([n])$ where
$$
\Phi ^{(m)}(t)=a_o^{q^m}+\sum \limits _{j=1}^ma_j^{q^m}\prod
\limits _{k=0}^{j-1}\left( t^{q^{m-k}}-\sum \limits
_{l=1}^k[1]^{q^{m-l}}\right)
$$
is a polynomial on $\K$ of a certain degree $N$ not depending on
$n$. Let $\theta _1,\ldots ,\theta _N$ be its roots. Then
$$
\Phi ^{(m)}([n])=a_m^{q^m}\prod \limits _{\nu =1}^N([n]-\theta _\nu
).
$$

As $n\to \infty $, $[n]\to -x$ in $\K$. We may assume that $\theta
_\nu \ne [n]$ for all $\nu$, if $n$ is large enough. If $\theta
_\nu \ne -x$ for all $\nu$, then for large $n$, say $n\ge n_0\ge
m$,
$$
\left| \Phi ^{(m)}([n])\right| \ge \mu >0.
$$
If $k\le N$ roots $\theta _\nu$ coincide with $-x$, then
$$
\left| \Phi ^{(m)}([n])\right| \ge \mu q^{-kq^n},\quad n\ge n_0.
$$

Combining the inequalities and taking the root we get
\begin{equation}
\left| \Phi _n\right| \ge \mu _1 q^{-\mu _2q^n},\quad n\ge n_0.
\end{equation}
where $\mu _1,\mu _2>0$. Now it follows from (12) and (13) that
the series (11) has (together with the series for $f$) a positive radius
of convergence.

Turning to the general equation (9) we note first of all that one
can apply an operator series $A(\tau )=\sum \limits _{k=0}^\infty
\alpha _k\tau ^k$ (even without assuming its convergence) to a
formal series (11), setting
$$
\tau ^ku(t)=\sum \limits _{n=0}^\infty u_n^{q^k}[n+1]^{q^{k-
1}}\ldots [n+k]\frac{t^{q^{n+k}}}{D_{n+k}},\quad k\ge 1,
$$
and
$$
A(\tau )u(t)=\sum \limits _{l=0}^\infty \frac{t^{q^l}}{D_l}\sum
\limits _{n+k=l}\alpha _ku_n^{q^k}[n+1]^{q^{k-1}}\ldots [n+k]
$$
where the factor $[n+1]^{q^{k-1}}\ldots [n+k]$ is omitted for
$k=0$.

Therefore the notion of a formal solution (11) makes sense for
the equation (9).

We will need the following elementary estimate.

\begin{lem}
Let $k\ge 2$ be a natural number, with a given partition
$k=i_1+\cdots +i_r$, where $i_1,\ldots ,i_r$ are positive
integers, $r\ge 1$. Then
$$
q^{i_1+\cdots +i_r}+q^{i_2+\cdots +i_r}+\cdots +q^{i_r}\le
q^{k+1}.
$$
\end{lem}

{\it Proof}. The assertion is obvious for $k=2$. Suppose it has
been proved for some $k$ and consider a partition
$$
k+1=i_1+\cdots +i_r.
$$
If $i_1>1$ then $k=(i_1-1)+i_2+\cdots +i_r$, so that
$$
q^{(i_1-1)+i_2+\cdots +i_r}+q^{i_2+\cdots +i_r}+\cdots +q^{i_r}\le
q^{k+1}
$$
whence
$$
q^{i_1+i_2+\cdots +i_r}+q^{i_2+\cdots +i_r}+\cdots +q^{i_r}\le
q^{k+2}.
$$

If $i_1=1$ then $k=i_2+\cdots +i_r$,
$$
q^{i_2+\cdots +i_r}+q^{i_3+\cdots +i_r}+\cdots +q^{i_r}\le
q^{k+1}
$$
and
$$
q^{i_1+\cdots +i_r}+q^{i_2+\cdots +i_r}+\cdots +q^{i_r}\le
2q^{k+1}\le q^{k+2}.\quad \blacksquare
$$

\medskip
Now we are ready to formulate our main result.

\begin{teo}
Let $u(t)$ be a formal solution (11) of the equation (9), where
the series for $A_j(\tau )z$, $z\in \K$, and $f(t)$, have
positive radii of convergence. Then the series (11) has a positive
radius of convergence.
\end{teo}

{\it Proof}. Applying (if necessary) the operator $\tau $ a
sufficient number of times to both sides of (9) we may assume
that
$$
A_j(\tau )=\sum \limits _{i=0}^\infty a_{ji}\tau ^{i+j},\quad
a_{ji}\in \K ,\ j=0,1,\ldots ,m,
$$
where $a_{j0}\ne 0$ at least for one value of $j$. Let us assume,
for example, that $a_{m0}\ne 0$ (otherwise the reasoning below
would need an obvious adjustment). Denote by $L$ the operator at
the left-hand side of (9), and by $L_0$ its ``principal part'',
$$
L_0u=\sum \limits _{j=0}^ma_{j0}\tau ^jd^ju
$$
(the model operator considered above; we will maintain the
notations introduced there). Note that $L_0$ is a linear
operator.

As we have seen,
$$
L_0\left( \frac{t^{q^n}}{D_n}\right) =\Phi _n\frac{t^{q^n}}{D_n},
\quad n\ge n_0,
$$
where $\Phi _n$ satisfies the inequality (13). This means that
$L_0$ is an automorphism of the vector space $X$ of formal series
$$
u=\sum \limits _{n=n_0}^\infty u_n\frac{t^{q^n}}{D_n},\quad
u_n\in \K,
$$
as well of its subspace $Y$ consisting of series with positive radii
of convergence.

Let us write the formal solution $u$ of the equation (9) as
$u=v+w$, where
$$
v=\sum \limits _{n=0}^{n_0-1} u_n\frac{t^{q^n}}{D_n}, \quad
w=\sum \limits _{n=n_0}^\infty u_n\frac{t^{q^n}}{D_n}.
$$
Then (9) takes the form
\begin{equation}
Lw=g,
\end{equation}
with $g=\sum g_n\frac{t^{q^n}}{D_n}\in Y$. In order to prove our
theorem, it is sufficient to verify that $w\in Y$.

For any $y\in X$ we can write
$$
Ly=(L_0-L_1)y=L_0(I-L_0^{-1}L_1)y
$$
where
\begin{equation}
L_1y=-\sum \limits _{j=0}^m\sum \limits _{i=1}^\infty a_{ji}\tau
^{i+j}d^jy.
\end{equation}
In particular, it is seen from (14) that
$$
(I-L_0^{-1}L_1)w=L_0^{-1}g,\quad L_0^{-1}g\in Y.
$$

Writing formally
$$
(I-L_0^{-1}L_1)^{-1}=\sum \limits _{k=0}^\infty \left( L_0^{-
1}L_1\right) ^k
$$
and noticing that $L_0^{-1}L_1:\ X\to \tau X$, we find that
\begin{equation}
w=\sum \limits _{k=0}^\infty \left( L_0^{-1}L_1\right) ^kh,
\end{equation}
where $h=L_0^{-1}g=\sum \limits _{n=n_0}^\infty
h_n\frac{t^{q^n}}{D_n}$, $h_n=\Phi _n^{-1}g_n$, and the series in
(16) converges in the natural non-Archimedean topology of the
space $X$.

A direct calculation shows that for any $\lambda \in \K$
$$
\left( L_0^{-1}L_1\right) \left( \lambda \frac{t^{q^n}}{D_n}\right) =-
\sum \limits _{i=1}^\infty \lambda ^{q^i}\Phi ^{-1}_{n+i}\Psi
_i^{(n)}\frac{t^{q^{n+i}}}{D_{n+i}}
$$
where
$$
\Psi _i^{(n)}=[n+1]^{q^{i-1}}[n+2]^{q^{i-2}}\ldots [n+i]\sum
_{j=0}^m[n-j+1]^{q^i}\ldots [n]^{q^i}a_{ji},
$$
and the coefficient at $a_{j0}$ in the last sum is assumed to
equal 1.

Proceeding by induction we get
\begin{multline*}
\left( L_0^{-1}L_1\right) ^r \left( \lambda \frac{t^{q^n}}{D_n}\right) =
(-1)^r\sum \limits _{i_1,\ldots ,i_r=1}^\infty \left( \Psi
_{i_1}^{(n)}\right) ^{q^{i_2+\cdots +i_r}}\left( \Psi
_{i_2}^{(n+i_1)}\right) ^{q^{i_3+\cdots +i_r}}\ldots \\
\times \left( \Psi _{i_r}^{(n+i_1+\cdots +i_{r-1})}\right)
\lambda ^{q^{i_1+\cdots +i_r}}\Phi _{n+i_1}^{-q^{i_2+\cdots +i_r}}
\Phi _{n+i_1+i_2}^{-q^{i_3+\cdots +i_r}}\ldots
\Phi ^{-1}_{n+i_1+\cdots +i_r}\frac{t^{q^{n+i_1+\cdots
+i_r}}}{D_{n+i_1+\cdots +i_r}},\quad r=1,2,\ldots .
\end{multline*}
Substituting this into (16) and changing the order of summation
we find an explicit formula for $w(t)$:
\begin{multline}
w(t)=\sum \limits _{l=n_0}^\infty \frac{t^{q^l}}{D_l}\sum \limits
_{\genfrac{}{}{0pt}{}{n+i_1+\cdots +i_r=l}{n\ge n_0,\, i_i,\ldots
,i_r\ge 1}}(-1)^rh_n^{q^{l-n}}\left( \Psi
_{i_1}^{(n)}\right) ^{q^{i_2+\cdots +i_r}}\left( \Psi
_{i_2}^{(n+i_1)}\right) ^{q^{i_3+\cdots +i_r}}\ldots \\
\times \left( \Psi _{i_r}^{(n+i_1+\cdots +i_{r-1})}\right)
\Phi _{n+i_1}^{-q^{i_2+\cdots +i_r}}
\Phi _{n+i_1+i_2}^{-q^{i_3+\cdots +i_r}}\ldots
\Phi ^{-1}_{n+i_1+\cdots +i_r}
\end{multline}

Observe that
$$
\left| \Psi _i^{(n)}\right| \le (q^{-1})^{q^{i-1}+q^{i-2}+\cdots
+1}\sup \limits _j|a_{ji}|,\quad |g_n|\le M_1^{q^n},\quad
|a_{ji}|\le M_2^{q^i},
$$
$M_1,M_2\ge 1$ (due to positivity of the corresponding radii of
convergence). We have
$$
\left| h_n^{q^{l-n}}\right| \le |\Phi _n|^{-q^{i_1+\cdots
+i_r}}M_1^{q^l},
$$
and by the above Lemma
\begin{multline*}
|\Phi _n|^{-q^{i_1+\cdots +i_r}}|\Phi _{n+i_1}|^{-q^{i_2+\cdots
+i_r}}\cdots |\Phi _{n+i_1+\cdots +i_r}|^{-1}\\
\le \mu _1^{q^{i_1+\cdots +i_r}+q^{i_2+\cdots +i_r}+\cdots
+q^{i_r}+1}q^{\mu _2\left( q^{n+i_1+\cdots +i_r}+q^{n+i_2+\cdots +i_r}+
\cdots +q^{n+i_r}+q^n\right) }\le
\mu _1^{q^{l-n+1}+1}q^{\mu_2q^n\left( q^{l-n+1}+1\right) }\\
\le q^{\mu _3q^{l+1}},\quad \mu_3>0.
\end{multline*}
The Lemma also yields
$$
\left| \Psi_{i_1}^{(n)}\right| ^{q^{i_2+\cdots +i_r}}\left| \Psi
_{i_2}^{(n+i_1)}\right| ^{q^{i_3+\cdots +i_r}}\ldots
\left| \Psi _{i_r}^{(n+i_1+\cdots +i_{r-1})}\right| \\
\le M_2^{q^{i_1+\cdots +i_r}+q^{i_2+\cdots +i_r}+\cdots
+q^{i_r}}\le M_2^{q^{l+1}}.
$$

Writing (17) as
$$
w(t)=\sum \limits _{l=n_0}^\infty w_l\frac{t^{q^l}}{D_l}
$$
we find that
$$
\limsup \limits_{l\to \infty}|w_l|^{q^{-l}}\le \limsup \limits_{l\to \infty}
\left( q^{\mu_3q^{l+1}}M_1^{q^l}M_2^{q^{l+1}}\right) ^{q{-
l}}<\infty ,
$$
which implies positivity of the radius of convergence. $\quad
\blacksquare $

\medskip
The function field analogue
$$
u(t)={}_2F_1(a,b;c;t),\quad a,b,c\in \mathbf Z,
$$
of the Gauss hypergeometric function, which was introduced by
Thakur [17], satisfies the equation [18]
$$
A_2(\tau )d^2u+A_1(\tau )du+A_0(\tau )u=0,
$$
where $A_2(\tau )=(1-\tau )\tau$, $A_1(\tau )=([-1]^q+[-b]+[-c])\tau
-[-c]$, $A_0(\tau )=-[-a][-b]$ (note that the element $[i]=x^{q^i}-x\in
\K$ is defined for any $i\in \mathbf Z$). The radii of convergence for
solutions of this equation are found explicitly in [18].

\section{ALGEBRAIC PROPERTIES}

In this section we consider the associative ring $\A$ of
``polynomial differential operators'', that is finite sums
\begin{equation}
a=\sum \limits _{i,j}\lambda_{ij}\tau ^id^j,\quad
\lambda_{ij}\in \K.
\end{equation}
Operations in $\A$ are defined in the natural way, with the use
of the commutation relations
$$
\tau \lambda =\lambda ^q\tau ,\quad d\lambda =\lambda ^{1/q}d\
(\lambda \in \K),\quad d\tau -\tau d=[1]^{1/q}.
$$

Note that a representation of an operator $a$ in the form (18) is
unique. Indeed, suppose that
$$
a=\sum \limits _{i=0}^m\sum \limits _{j=0}^n\lambda_{ij}\tau
^id^j=0.
$$
Let $\psi _l(t)=\frac{t^{q^l}}{D_l}$. Using (1), (2), we find
that
$$
0=a(\psi _0)=\sum \limits _{i=0}^m\lambda _{i0}\tau ^i\psi
_0=\lambda _{00}\psi _0+\sum \limits _{i=1}^m\lambda
_{i0}[1]^{q^{i-1}}[2]^{q^{i-2}}\ldots [i]\psi _i
$$
whence $\lambda _{i0}=0$. Then we proceed by induction; if
$\lambda _{ij}=0$ for $j\le \nu<n$, then
$$
a(\psi _{\nu +1})=\sum \limits _{i=0}^m\sum \limits _{j=\nu +1}^n
\lambda_{ij}\tau ^id^j\psi _{\nu +1}=\sum \limits _{i=0}^m\lambda
_{i,\nu +1}\tau ^i\psi_0,
$$
so that $\lambda _{i,\nu +1}=0$ ($i=0,1,\ldots ,m$) as before.

Some algebraic properties of the ring $\A$ are collected in the
following theorem.

\begin{teo}
(i) The centre of the ring $\A$ coincides with $\F$.

(ii) The ring $\A$ possesses no non-trivial two-sided ideals
stable with respect to the mapping
$$
P\left( \sum \limits _{i,j}\lambda_{ij}\tau ^id^j\right) =
\sum \limits _{i,j}\lambda_{ij}^q\tau ^id^j.
$$

(iii) The ring $\A$ is Noetherian.

(iv) $\A$ is an Ore domain, that is $\A$ has no zero-divisors and
$\A a\cap \A b\ne \{ 0\}$, $a\A \cap b\A \ne \{ 0\}$ for all
pairs of non-zero elements $a,b\in \A$.
\end{teo}

{\it Proof}. (i) It is easily proved (by induction) that
\begin{equation}
[d,\tau ^i]=[i]^{1/q}\tau ^{i-1},\quad [d^j,\tau ]=[j]^{q^{-
j}}d^{j-1}
\end{equation}
for any natural numbers $i,j$.

Suppose that $a=\sum \limits _{i=0}^m\sum \limits _{j=0}^n\lambda_{ij}\tau
^id^j$ belongs to the centre of $\A$. Then $[\tau ,a]=[d,a]=0$.
Repeatedly using (19), we find that
$$
0=[\tau ,a]=\sum \limits _{i=0}^m\sum \limits _{j=0}^n
\left( \lambda_{ij}^q-\lambda _{ij}\right) \tau^{i+1}d^j-\sum \limits
_{j=1}^n\lambda _{0j}[j]^{q^{-j}}d^{j-1}-\sum \limits _{i=0}^{m-1}
\sum \limits _{j=0}^{n-1}\lambda_{i+1,j+1}[j+1]^{q^{i-j}}
\tau^{i+1}d^j
$$
whence
\begin{equation}
\lambda _{0j}=0,\quad j=1,\ldots ,n,
\end{equation}
\begin{equation}
\lambda _{ij}^q-\lambda _{ij}-[j+1]^{q^{i-j}}\lambda _{i+1,j+1}=0,
\quad i=0,1,\ldots ,m-1;\ j=0,1,\ldots ,n-1.
\end{equation}

It follows from (20), (21) that $\lambda _{ij}=0$ for $i<j$.
Next,
\begin{multline*}
0=[d,a]=\sum \limits _{i=1}^m\sum \limits _{j=0}^i
\left( \lambda_{ij}^{1/q}-\lambda _{ij}\right) \tau^id^{j+1}
-\sum \limits _{i=0}^{m-1}\sum \limits _{j=0}^i\lambda_
{i+1,j+1}^{1/q}[i+1]^{1/q}\tau^id^{j+1}\\
+\left( \lambda _{00}^{1/q}-\lambda _{00}\right) d-\sum \limits
_{i=1}^m\lambda _{i0}^{1/q}[i]^{1/q}\tau ^{i-1},
\end{multline*}
so that
\begin{equation}
\lambda _{i0}=0,\quad i=1,\ldots ,n;
\end{equation}
\begin{equation}
\lambda _{ij}^{1/q}-\lambda _{ij}-\lambda_{i+1,j+1}^{1/q}[i+1]^{1/q}=0,
\quad i=1,\ldots ,n-1;\ j=0,1,\ldots ,i;
\end{equation}
\begin{equation}
\lambda _{00}^{1/q}-\lambda _{00}-\lambda_{11}^{1/q}[1]^{1/q}=0.
\end{equation}

From (22), (23) we get $\lambda _{ij}=0$ for $i>j$. Raising (24)
to the power $q$ we can compare the resulting equality with (21)
(with $i=j=0$). Then we find that $\lambda _{11}=0$, and by
virtue of (21) also $\lambda _{ii}=0$, $i=2,\ldots ,m$. Finally,
it follows from (24) that $\lambda _{00}^q=\lambda _{00}$, so
that $a=\lambda _{00}\in \F$.

(ii) Let $D$ be a two-sided ideal in $\A$, $PD\subset D$,
containing a non-zero element
$$
a=\sum \limits _{i=0}^m\sum \limits _{j=0}^n\lambda_{ij}\tau
^id^j.
$$
Then $D$ contains the element $a_1=P(a\tau )-\tau a$. It follows
as above that
$$
a_1=\sum \limits _{i=0}^m\sum \limits _{j=1}^n\lambda_{ij}^q[j]^{i-j+1}\tau
^id^{j-1}.
$$

It is clear that either $\lambda _{ij}=0$ for $j\ge 1$, or
$a_1\ne 0$, and the maximal degree of $d$ in $a_1$ is smaller by
1 than the one in $a$. Repeating the procedure (if necessary) we
obtain a non-zero element of $D$ of the form
$$
b=\sum \limits _{i=0}^m\mu _i\tau ^i.
$$

If not all the coefficients $\mu _i$, $i\ge 1$, are equal to
zero, we find a non-zero $b_1\in D$, $b_1=P(db)-bd$,
$$
b_1=\sum \limits _{i=1}^m\mu _i[i]\tau ^{i-1}.
$$
After an appropriate repetition we obtain that $D$ contains a
non-zero constant, so that $D=\A$.

(iii) Let
$$
A_\nu =\left\{ \sum \limits _{i=0}^m\sum \limits _{j=0}^n\lambda_{ij}\tau
^id^j\in \A\ :\ m+n\le \nu \right\} .
$$
The sequence $\{ A_\nu \}$ of $\K$-vector spaces is increasing,
and we can define a graded ring
$$
\G =A_0\oplus A(1)\oplus A(2)\oplus \ldots
$$
where $A(\nu )=A_\nu /A_{\nu -1}$, $\nu \ge 1$. The
multiplication in $\G$ is defined as follows. If $f\in A(\nu )$,
$g\in A(k)$, $\varphi \in A_\nu$ and $\psi \in A_k$ are arbitrary
representatives of $f$ and $g$ respectively, then $\varphi \psi
\in A_{\nu +k}$, and we define $fg$ as the class of $\varphi
\psi$ in $A(\nu +k)$. In can be checked easily that the
multiplication is well-defined.

The ring $\G$ is generated by the classes $\overline{\tau
},\overline{d}\in A(1)$ of the elements $\tau ,d\in A_1$, and
constants from $\K$, with the commutation relations
$$
\overline{\tau }\overline{d}-\overline{d}\overline{\tau }=0,
\quad \overline{\tau }c=c^q\overline{\tau }, \quad
\overline{d}c=c^{1/q}\overline{d}\ \ (c\in \K ).
$$
It follows from the generalization of the Hilbert basis theorem
given in \cite{Rt} that $\G$ is a Noetherian ring.

Let $\LL$ be a left ideal in $\A$. We have to prove that $\LL$ is
finitely generated as an $\A$-submodule. Let $\Gamma _\nu =A_\nu
\cap \LL$. Then $\{ \Gamma _\nu \}$ is a filtration in $\LL$. As
above, we can construct a graded ring $\GL$, which is a left
ideal in $\G$. Since $\G$ is Noetherian, $\GL$ has a finite system
of generators $\sigma _1,\ldots ,\sigma _n$, and we can write
finite decompositions
$$
\sigma _j=\sum \limits _k\sigma _j(k),\quad \sigma _j(k)\in
\Gamma _k/\Gamma _{k-1}.
$$
Let $\mu _1,\ldots ,\mu _N$ be the set of all non-zero elements
$\sigma _j(k)$.

Denote by $\gamma $ the canonical imbedding $\Gamma _k\to
\Gamma _k/\Gamma _{k-1}=\Gamma (k)$ extended to the mapping $\LL
\to \GL$. Choose $m_i\in \LL$ in such a way that $\gamma
(m_i)=\mu _i$. Then $m_i$ are generators of $\LL$, that is
\begin{equation}
\Gamma _k\subset \sum \limits _{i=1}^N\A m_i\quad \mbox{for each
}k.
\end{equation}

Indeed, (25) is obvious for $k=0$. Suppose that
\begin{equation}
\Gamma _{k-1}\subset \sum \limits _{i=1}^N\A m_i\quad \mbox{for each
}k.
\end{equation}
Consider an element $l\in \Gamma _k\setminus \Gamma _{k-1}$. We
have $\gamma (l)\in \Gamma (k)$,
$$
\gamma (l)=\sum \limits _{j=1}^nc_j\sigma _j, \quad c_j\in \G .
$$
Writing each $c_j$ as a sum of homogeneous components
$$
c_j=\sum \limits _\nu c_j(\nu ),\quad c_j(\nu )\in A(\nu ),
$$
and taking
into account that $A(\nu )\Gamma (k)\subset \Gamma (k+\nu )$ for
any $k,\nu $, we find that
$$
\gamma (l)=\sum \limits _{j+\nu =k}c_j(\nu )\sigma _j(k)=
\sum \limits _{j=0}^kc_j(k-j)\sigma _j(k).
$$
Choosing $C_j\in A_{k-j}$ in such a way that $c_j(k-j)$ is a
class of $C_j$ in $A(k-j)$, we obtain the inclusion
$$
l-\sum \limits _{j=0}^kC_jm'_j\in \Gamma _{k-1}
$$
where $\{ m'_j\}$ is a subset of $\{ m_i\} _{i=1}^N$. Together
with (26) this implies (25).

We have proved that $\A$ is left Noetherian. The proof of the
right Noetherianness is similar.

(iv) Let $ab=0$ for
$$
a=\sum \limits _{i=0}^{m_1}\sum \limits _{j=0}^{n_1}\lambda_{ij}\tau
^id^j,\quad b=\sum \limits _{k=0}^{m_2}\sum \limits _{l=0}^{n_2}
\mu_{kl}\tau ^kd^l,
$$
and $a\ne 0$, $b\ne 0$, that is
\begin{equation}
\sum \limits _{i=0}^{m_1}\lambda _{in_1}\tau ^i\ne 0,\quad
\sum \limits _{k=0}^{m_2}\mu _{kn_2}\tau ^k\ne 0.
\end{equation}

It follows from (19) that
$$
d^j\tau ^k=\tau ^kd^j+O\left( d^{j-1}\right)
$$
where $O\left( d^{j-1}\right)$ means a polynomial in the variable $d$,
of a degree $\le j-1$, with coefficients from the composition ring
$\K \{ \tau \}$ of polynomials in the operator $\tau $.

Therefore the coefficient at $d^{n_1+n_2}$ in the expression for the operator
$ab$  equals
$$
\sum \limits _{i,k}\lambda _{in_1}\mu _{kn_2}^{q^{i-n_1}}\tau ^{i+k}=
P_1(P_2(\tau ))
$$
where
$$
P_1(\tau )=\sum \limits _{i=0}^{m_1}\lambda _{in_1}\tau ^i,\quad
P_2(\tau )=\sum \limits _{k=0}^{m_2}\mu _{kn_2}^{q^{-n_1}}\tau ^k,
$$
which contradicts (27), since the ring $\K \{ \tau \}$ has no
zero-divisors \cite{G2}.

Now (iv) follows from (iii) (see Sect. 4.5 in \cite{Bok}).
However we will give also an elementary direct proof (which does
not use the Hilbert basis theorem or its generalizations).

Let us prove the left Ore condition (the proof of the right
condition is simpler and does not differ from the one in
characteristic zero, see \cite{Bj}). Thus let $a,b\ne 0$; we will
prove that $a\A \cap b\A \ne \{ 0\}$.

As above, we will use the filtration $\{ A_\nu \}$ in $\A$.
Let $\nu _1$ be such a number that $a,b\in A_{\nu _1}$. Then
$aA_\nu \subset A_{\nu +\nu _1}$, $bA_\nu \subset A_{\nu +\nu
_1}$. Suppose that $a\A \cap b\A =\{ 0\}$. Then in particular
\begin{equation}
aA_\nu \cap bA_\nu =\{ 0\}
\end{equation}

Let us prove that the set $\{ a\tau ^id^j,\ i+j\le \nu \}$ forms
a basis of $aA_\nu$, that is its elements are linearly
independent. This assertion is evident for $\nu =0$. Suppose that
it holds for all $\nu \le N-1$, and consider the case $\nu =N$.
Let
\begin{equation}
\sum \limits _{i+j\le N}c_{ij}a\tau ^id^j=0,\quad c_{ij}\in \K.
\end{equation}

We may write
$$
a=\sum \limits _{k+l\le \kappa }\lambda _{kl}\tau ^kd^l,\quad
\lambda _{kl}\in \K,\ \kappa \le \nu _1,
$$
where $\lambda _{kl}\ne 0$ at least for one couple $(k,l)$ with
$k+l=\kappa$. Since
$$
d^l\tau^i\equiv \tau ^id^l\pmod{A_{i+l-1}},
$$
we find that
$$
\sum \limits _{i+j=N}c_{ij}\left( \sum \limits _{k+l=
\kappa }\lambda _{kl}\tau ^kd^l\right) \tau ^id^j\equiv
\sum \limits _{i+j=N}\sum \limits _{k+l=\kappa }c_{ij}\lambda _{kl}
\tau ^{i+k}d^{j+l}\pmod{A_{i+k+j+l-1}}.
$$
By (29), this means that
$$
\sum \limits _{i+j=N}\sum \limits _{k+l=\kappa }c_{ij}\lambda _{kl}
\tau ^{i+k}d^{j+l}=0
$$
whence
$$
0=\sum \limits _{i=0}^Nc_{i,N-i}\sum \limits _{k=0}^\kappa
\lambda _{k,\kappa -k}\tau ^{i+k}d^{N+\kappa -(i+k)}=\sum \limits
_{m=0}^{N+\kappa }\left( \sum \limits _{i+k=m}c_{i,N-i}
\lambda _{k,\kappa -k}\right) \tau ^md^{N+\kappa -m},
$$
so that
$$
\sum \limits _{i+k=m}c_{i,N-i}\lambda _{k,\kappa -k}=0,\quad
m=0,1,\ldots ,N+\kappa .
$$
The expression in the left-hand side coincides with the $m$-th
coefficient of the product of two polynomials. Therefore
$c_{ij}=0$ for $i+j=N$, the summation in (29) is actually
performed for $i+j\le N-1$, and by the induction assumption
$c_{ij}=0$ for all $i,j$.

Now $\dim (aA_\nu )=\dim (bA_\nu)=\dim A_\nu $, and it follows
from (28) that
\begin{equation}
\dim (A_{\nu +\nu _1})\ge \dim (aA_\nu \oplus bA_\nu)=2\dim
(A_\nu ).
\end{equation}
Note that
$$
\dim A_\nu =\mbox{card }\{ (i,j):\ i+j\le \nu \} =\frac{(\nu
+1)(\nu +2)}2,
$$
and we see that
$$
\frac{\dim (A_{\nu +\nu _1})}{\dim A_\nu }\longrightarrow 1\ \
\mbox{as}\ \ \nu \to \infty ,
$$
which contradicts (30). $\qquad \blacksquare$

\end{document}